\documentclass{article}%
\usepackage{amsfonts}
\usepackage{amsmath}%
\usepackage{amssymb}%
\usepackage{graphicx}

\begin{document}

\title{Chromogeometry and relativistic conics}
\date{}
\author{N J Wildberger\\School of Mathematics and Statistics\\UNSW Sydney 2052 Australia}
\maketitle

This paper shows how a recent reformulation of the basics of classical
geometry and trigonometry reveals a three-fold symmetry between Euclidean and
non-Euclidean (relativistic) planar geometries. We apply this
\textit{chromogeometry} to look at conics in a new light.

\section*{Pythagoras, area and quadrance}

To measure a line segment in the plane, the ancient Greeks measured the
\textit{area of a square constructed on it}. Algebraically, the parallelogram
formed by a vector $v=\overrightarrow{A_{1}A_{2}}=\left(  a,b\right)  $ and
its perpendicular $B\left(  v\right)  =\left(  -b,a\right)  $ has area%
\begin{equation}
Q=\det%
\begin{pmatrix}
a & b\\
-b & a
\end{pmatrix}
=a^{2}+b^{2}.
\end{equation}
The Greeks referred to building squares as `quadrature', and so we say that
$Q$ is the \textbf{quadrance} of the vector $v,$ or the quadrance $Q\left(
A_{1},A_{2}\right)  $ between $A_{1}$ and $A_{2}.$ This notion makes sense
over any field.%
\begin{figure}
[h]
\begin{center}
\includegraphics[
height=3.8785cm,
width=4.488cm
]%
{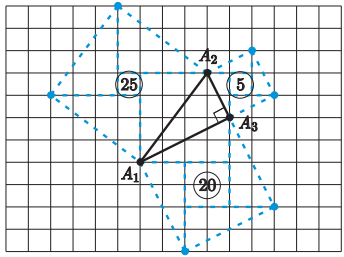}%
\caption{Pythagoras' theorem: $5+20=25$}%
\label{BluePythagoras}%
\end{center}
\end{figure}

If $Q_{1}=Q\left(  A_{2},A_{3}\right)  $, $Q_{2}=Q\left(  A_{1},A_{3}\right)
$ and $Q_{3}=Q\left(  A_{1},A_{2}\right)  $ are the quadrances of a triangle
$\overline{A_{1}A_{2}A_{3}},$ then Pythagoras' theorem and its converse can
together be stated as: $A_{1}A_{3}$\textit{\ is perpendicular to }$A_{2}A_{3}
$\textit{\ precisely when}%
\[
Q_{1}+Q_{2}=Q_{3}.
\]
Figure \ref{BluePythagoras} shows an example where $Q_{1}=5,$ $Q_{2}=20$ and
$Q_{3}=25.$ As indicated for the large square, these areas may also be
calculated by subdivision and (translational) rearrangement, followed by
counting cells.

\begin{figure}
[h]
\begin{center}
\includegraphics[
height=3.7752cm,
width=4.9424cm
]%
{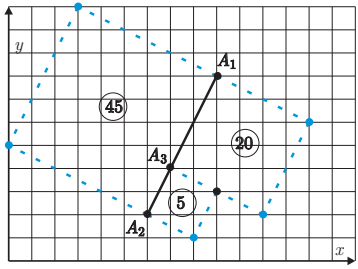}%
\caption{Triple quad formula: $\left(  5+20+45\right)  ^{2}=2\left(
5^{2}+20^{2}+45^{2}\right)  $ }%
\label{BlueTriple}%
\end{center}
\end{figure}

There is a sister theorem---the \textit{Triple quad formula}---that Euclid did
not know, but which is fundamental for \textit{rational trigonometry},
introduced in \cite{Wild1}: $A_{1}A_{3}$\textit{\ is parallel to }$A_{2}A_{3}
$\textit{\ precisely when}%
\[
\left(  Q_{1}+Q_{2}+Q_{3}\right)  ^{2}=2\left(  Q_{1}^{2}+Q_{2}^{2}+Q_{3}%
^{2}\right)  .
\]
Figure \ref{BlueTriple} shows an example where $Q_{1}=5,$ $Q_{2}=20$ and
$Q_{3}=45.$%

In terms of side lengths $d_{1}=\sqrt{Q_{1}},d_{2}=\sqrt{Q_{2}}$ and
$d_{3}=\sqrt{Q_{3}},$ and the semi-perimeter $s=\left(  d_{1}+d_{2}%
+d_{3}\right)  /2,$ observe that%
\begin{align*}
&  \left(  Q_{1}+Q_{2}+Q_{3}\right)  ^{2}-2\left(  Q_{1}^{2}+Q_{2}^{2}%
+Q_{3}^{2}\right) \\
&  =4Q_{1}Q_{2}-\left(  Q_{1}+Q_{2}-Q_{3}\right)  ^{2}=4d_{1}^{2}d_{2}%
^{2}-\left(  d_{1}^{2}+d_{2}^{2}-d_{3}^{2}\right)  ^{2}\\
&  =\left(  2d_{1}d_{2}-\left(  d_{1}^{2}+d_{2}^{2}-d_{3}^{2}\right)  \right)
\left(  2d_{1}d_{2}+\left(  d_{1}^{2}+d_{2}^{2}-d_{3}^{2}\right)  \right) \\
&  =\left(  d_{3}^{2}-\left(  d_{1}-d_{2}\right)  ^{2}\right)  \left(  \left(
d_{1}+d_{2}\right)  ^{2}-d_{3}^{2}\right) \\
&  =\left(  d_{3}-d_{1}+d_{2}\right)  \allowbreak\left(  d_{3}+d_{1}%
-d_{2}\right)  \left(  d_{1}+d_{2}-d_{3}\right)  \left(  d_{1}+d_{2}%
+d_{3}\right) \\
&  =16\left(  s-d_{1}\right)  \left(  s-d_{2}\right)  \left(  s-d_{3}\right)
s.
\end{align*}
Thus Heron's formula in the usual form%
\[
\mathrm{area}=\sqrt{s\left(  s-d_{1}\right)  \left(  s-d_{2}\right)  \left(
s-d_{3}\right)  }%
\]
may be restated in terms of quadrances as%
\[
16~\mathrm{area}^{2}=\left(  Q_{1}+Q_{2}+Q_{3}\right)  ^{2}-2\left(  Q_{1}%
^{2}+Q_{2}^{2}+Q_{3}^{2}\right)  \equiv A\left(  Q_{1},Q_{2},Q_{3}\right)  .
\]
This more fundamental formulation deserves to be called \textbf{Archimedes'
theorem}, since Arab sources indicate that Archimedes knew Heron's formula.
The Triple quad formula is the special case of Archimedes' theorem when the
area is zero. The function $A\left(  Q_{1},Q_{2},Q_{3}\right)  $ will be
called \textbf{Archimedes' function}.

\begin{figure}
[h]
\begin{center}
\includegraphics[
height=3.6416cm,
width=5.1491cm
]%
{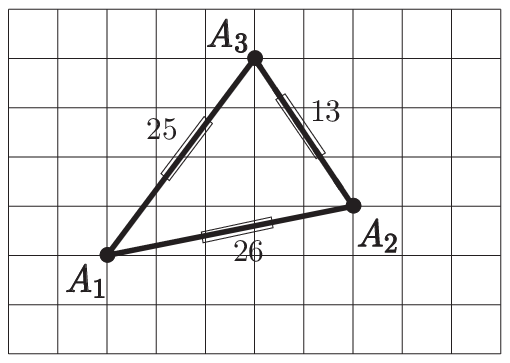}%
\caption{$16~\mathrm{area}^{2}=\left(  13+25+26\right)  ^{2}-2\left(
13^{2}+25^{2}+26^{2}\right)  =\allowbreak1156$}%
\label{Arch}%
\end{center}
\end{figure}

In Figure \ref{Arch} the quadrances are $Q_{1}=13,$ $Q_{2}=25$ and $Q_{3}=26,
$ so $16~\mathrm{area}^{2}=\allowbreak1156$, giving $\mathrm{area}^{2}=289/4 $
and $\mathrm{area}=17/2.$
Irrational side lengths are not needed to determine the area of a rational
triangle, and in any case when we move to more general geometries, we have no
choice but to give up on distance and angle.

\section*{Blue, red and green geometries}

Euclidean geometry will here be called \textbf{blue geometry}. We now
introduce two relativistic geometries, called\textit{\ red} and \textit{green}%
, which arise from Einstein's theory of relativity. These rest on alternate
notions of perpendicularity, but they share the same underlying affine concept
of area as blue geometry, and indeed the same laws of rational trigonometry,
as will be explained shortly.%
\begin{figure}
[h]
\begin{center}
\includegraphics[
height=1.5333in,
width=3.8623in
]%
{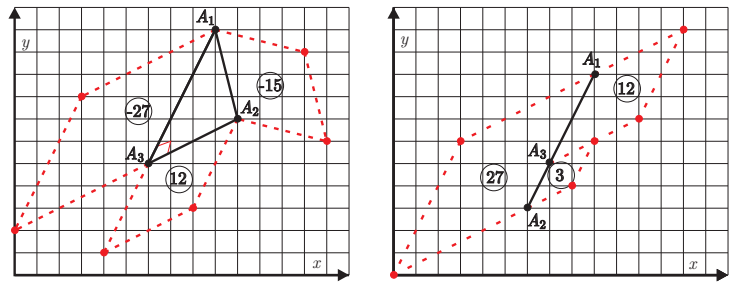}%
\caption{Red Pythagoras' theorem and Triple quad formula}%
\label{RedPythTriple}%
\end{center}
\end{figure}

Define the vector $v=\overrightarrow{A_{1}A_{2}}=\left(  a,b\right)  $ to be
\textbf{red perpendicular} to $R\left(  v\right)  =\left(  b,a\right)  .$ This
mapping is easily visualized: it corresponds to Euclidean reflection in a line
of slope $1$ or $-1.$ A \textbf{red square} is then a parallelogram with sides
$v$ and $R\left(  v\right)  $, and hence (signed) area
\begin{equation}
Q^{\left(  r\right)  }=\det%
\begin{pmatrix}
a & b\\
b & a
\end{pmatrix}
=a^{2}-b^{2}%
\end{equation}
which we call the \textbf{red quadrance} between $A_{1}$ and $A_{2}.$ Figure
\ref{RedPythTriple} illustrates that both Pythagoras' theorem and the Triple
quad formula hold also using red quadrances and red perpendicularity, where
the areas of the red squares can be computed as before by subdivisions,
(translational) rearrangement and counting cells---or by applying the
algebraic formula for the red quadrance.

In a similar fashion the vector $v=\overrightarrow{A_{1}A_{2}}=\left(
a,b\right)  $ is \textbf{green perpendicular} to $G\left(  v\right)  =\left(
-a,b\right)  $. This corresponds to Euclidean reflection in a vertical or
horizontal line. A \textbf{green square} is a parallelogram with sides $v$ and
$R\left(  v\right)  $, and hence (signed) area
\begin{equation}
Q^{\left(  g\right)  }=\det%
\begin{pmatrix}
a & b\\
-a & b
\end{pmatrix}
=2ab
\end{equation}
which we call the \textbf{green quadrance} between $A_{1}$ and $A_{2}.$%
\begin{figure}
[h]
\begin{center}
\includegraphics[
height=3.5727cm,
width=10.9292cm
]%
{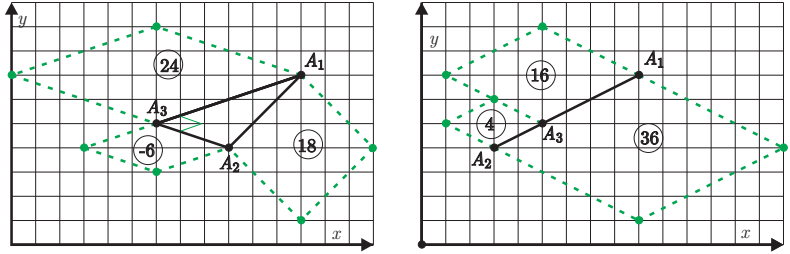}%
\caption{Green\ Pythagoras' theorem and Triple quad formula}%
\label{GreenPythTriple}%
\end{center}
\end{figure}

Figure \ref{GreenPythTriple} shows Pythagoras' theorem and the Triple quad
formula in the green context. This version of relativistic geometry
corresponds to a basis of null vectors in red geometry.

All three geometries can be defined over a general field, not of
characteristic two.

\section*{Spreads and rational trigonometry}

The three quadratic forms%
\[%
\begin{tabular}
[c]{lllll}%
$Q^{\left(  b\right)  }\left(  a,b\right)  =a^{2}+b^{2}$ &  & $Q^{\left(
r\right)  }\left(  a,b\right)  =a^{2}-b^{2}$ &  & $Q^{\left(  g\right)
}\left(  a,b\right)  =2ab$%
\end{tabular}
\]
have corresponding dot products%
\begin{align*}
\left(  a_{1},b_{1}\right)  \cdot_{b}\left(  a_{2},b_{2}\right)   & \equiv
a_{1}a_{2}+b_{1}b_{2}\\
\left(  a_{1},b_{1}\right)  \cdot_{r}\left(  a_{2},b_{2}\right)   & \equiv
a_{1}a_{2}-b_{1}b_{2}\\
\left(  a_{1},b_{1}\right)  \cdot_{g}\left(  a_{2},b_{2}\right)   & \equiv
a_{1}b_{2}+a_{2}b_{1}.
\end{align*}
Together with $a_{1}b_{2}-a_{2}b_{1}=0$ describing parallel vectors, these are
the four simplest bilinear expressions in the four variables.

In rational trigonometry, one wants to work over general fields, so the notion
of angle is not available, but it is important to realize that the dot product
is not necessarily the best replacement. Instead we introduce the related
notion of \textit{spread }between two lines (not between rays), which in the
blue framework is the square of the sine of the angle between the lines (there
are actually many such angles, but the square of the sine is the same for all).

The \textbf{blue, red }and \textbf{green spreads} between lines $l_{1}$ and
$l_{2}$ with equations
\[%
\begin{tabular}
[c]{lllll}%
$a_{1}x+b_{1}y+c_{1}=0$ &  & \textrm{and} &  & $a_{2}x+b_{2}y+c_{2}=0$%
\end{tabular}
\]
are respectively the numbers
\begin{align*}
s^{\left(  b\right)  }\left(  l_{1},l_{2}\right)   & =\frac{\left(  a_{1}%
b_{2}-a_{2}b_{1}\right)  ^{2}}{\left(  a_{1}^{2}+b_{1}^{2}\right)  \left(
a_{2}^{2}+b_{2}^{2}\right)  }\\
s^{\left(  r\right)  }\left(  l_{1},l_{2}\right)   & =-\frac{\left(
a_{1}b_{2}-a_{2}b_{1}\right)  ^{2}}{\left(  a_{1}^{2}-b_{1}^{2}\right)
\left(  a_{2}^{2}-b_{2}^{2}\right)  }\\
s^{\left(  g\right)  }\left(  l_{1},l_{2}\right)   & =-\frac{\left(
a_{1}b_{2}-a_{2}b_{1}\right)  ^{2}}{4a_{1}b_{1}a_{2}b_{2}}.
\end{align*}
These quantities are undefined when the denominators are zero. The negative
signs in front of $s^{\left(  r\right)  }$ and $s^{\left(  g\right)  }$ insure
that, for each of the colours, the \textit{spread at any of the three vertices
of a right triangle (one with two sides perpendicular) is the quotient of the
opposite quadrance by the hypotenuse quadrance}. See \cite{Wild2} for a proof
of this, and other facts about rational trigonometry, in a wider context.

In Figure \ref{BluePythagoras} the spreads at $A_{1}$ and $A_{2}$ are $1/5$
and $4/5$ respectively, in the left diagram of Figure \ref{RedPythTriple} the
spreads at $A_{1}$ and $A_{2}$ are $-4/5$ and $9/5$ respectively, and in the
left diagram of Figure \ref{GreenPythTriple} the spreads at $A_{1}$ and
$A_{2}$ are $-1/3$ and $4/3$ respectively. In each case the spread at the
right vertex is $1$, and the other two spreads sum to $1.$%
\begin{figure}
[h]
\begin{center}
\includegraphics[
height=3.689cm,
width=10.744cm
]%
{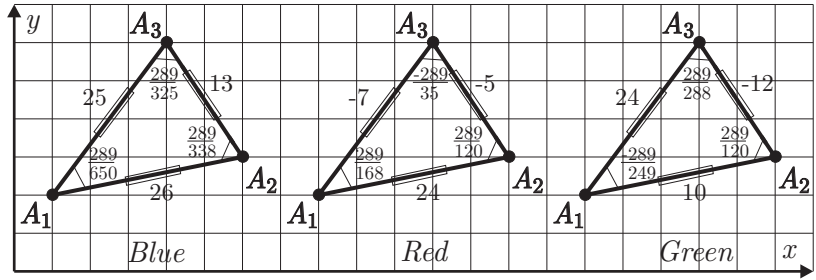}%
\caption{Blue, red and green quadrances and spreads}%
\label{TriangleExample}%
\end{center}
\end{figure}

Figure \ref{TriangleExample} allows you to compare the various quadrances and
spreads of a fixed triangle in each of the three geometries. Note the common
numerators of the spreads, arising because $A\left(  Q_{1},Q_{2},Q_{3}\right)
=\pm4\times289$ is up to sign the same in each geometry, with a plus sign in
the blue situation and a negative sign in the red and green ones.

Aside from Pythagoras' theorem and the Triple quad formula, the main laws of
rational trigonometry are: for a triangle with quadrances $Q_{1},Q_{2}$ and
$Q_{3},$ and spreads $s_{1},s_{2}$ and $s_{3}:$
\[
\frac{s_{1}}{Q_{1}}=\frac{s_{2}}{Q_{2}}=\frac{s_{3}}{Q_{3}}%
~~\text{\textrm{(Spread law)}}%
\]%
\[
\left(  Q_{1}+Q_{2}-Q_{3}\right)  ^{2}=4Q_{1}Q_{2}\left(  1-s_{3}\right)
~~\text{\textrm{(Cross law)}}%
\]%
\[
\left(  s_{1}+s_{2}+s_{3}\right)  ^{2}=2\left(  s_{1}^{2}+s_{2}^{2}+s_{3}%
^{2}\right)  +4s_{1}s_{2}s_{3}~~\text{\textrm{(Triple spread formula).}}%
\]
As shown in \cite{Wild1}, these laws are derived using only Pythagoras'
theorem and the Triple quad formula. Since these latter two results hold in
all three geometries, the Spread law, Cross law and Triple spread formula also
hold in all three geometries.

For any points $A_{1}$ and $A_{2}$ the square of $Q^{\left(  b\right)
}\left(  A_{1},A_{2}\right)  $ is the sum of the squares of $Q^{\left(
r\right)  }\left(  A_{1},A_{2}\right)  $ and $Q^{\left(  g\right)  }\left(
A_{1},A_{2}\right)  $, and for any lines $l_{1}$ and $l_{2}$
\[
\frac{1}{s^{\left(  b\right)  }\left(  l_{1},l_{2}\right)  }+\frac
{1}{s^{\left(  r\right)  }\left(  l_{1},l_{2}\right)  }+\frac{1}{s^{\left(
g\right)  }\left(  l_{1},l_{2}\right)  }=2.
\]
The first statement follows from the Pythagorean triple identity
\[
\left(  x^{2}+y^{2}\right)  ^{2}=\left(  x^{2}-y^{2}\right)  ^{2}+\left(
2xy\right)  ^{2}%
\]
while the latter follows from the identity%
\[
\left(  a_{1}^{2}+b_{1}^{2}\right)  \left(  a_{2}^{2}+b_{2}^{2}\right)
-\left(  a_{1}^{2}-b_{1}^{2}\right)  \left(  a_{2}^{2}-b_{2}^{2}\right)
-4a_{1}b_{1}a_{2}b_{2}=\allowbreak2\left(  a_{1}b_{2}-a_{2}b_{1}\right)  ^{2}.
\]
So in Figure \ref{TriangleExample} there are three linked (signed) Pythagorean
triples, namely $\left(  13,-5,-12\right)  $, $\left(  25,-7,24\right)  $ and
$\left(  26,24,10\right)  $, and three triples of harmonically related
spreads.%
\begin{figure}
[h]
\begin{center}
\includegraphics[
height=3.6115cm,
width=5.147cm
]%
{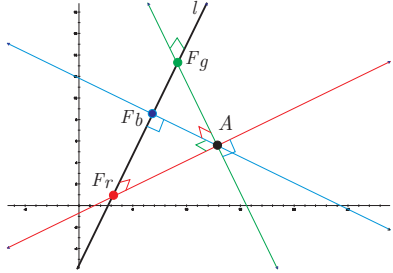}%
\caption{Three altitudes from a point to a line}%
\label{Altitudes}%
\end{center}
\end{figure}

Figure \ref{Altitudes} shows the three coloured altitudes from a point $A$ to
a line $l,$ and the feet of those altitudes. Note that the blue and red
altitudes are green perpendicular and similarly for the other colours. The
three triangles formed by the four points are each \textbf{triple right
triangles}, containing each a blue, red and green right vertex.

Most theorems of planar Euclidean geometry have universal versions, valid in
each of the three geometries. \textit{This is a large claim that deserves
further investigation}. In the red and green geometries, circles and rotations
become rectangular hyperbolas and Lorentz boosts. There are no equilateral
triangles in the red and green geometries, so results like Napoleon's theorem
or Morley's theorem will not have (obvious) analogs.

To see some chromogeometry in action, let's have a look at conics\textit{\ }in
this more general framework.

\section*{The ellipse as a grammola}

In the real number plane, one usually defines an \textit{ellipse} as the locus
of a point $X$ whose ratio of distance from a fixed point (focus) to distance
from a fixed line (directrix) is constant and less than one, and hyperbolas
and parabolas similarly with eccentricities greater than one and equal to one.
By squaring this condition, we can discuss the locus of a point whose ratio of
\textit{quadrance} from a fixed point to a fixed line is constant. By
quadrance from a point $X$ to a line $l$ we mean the obvious: construct the
altitude line $n$ from $X$ to $l,$ find its foot $F$ and measure $Q\left(
X,F\right)  $. Let's call such a locus a \textbf{conic section}. Over a
general field we cannot distinguish `ellipses' from `hyperbolas', although
parabolas are always well defined.
\begin{figure}
[h]
\begin{center}
\includegraphics[
height=5.1814cm,
width=11.0606cm
]%
{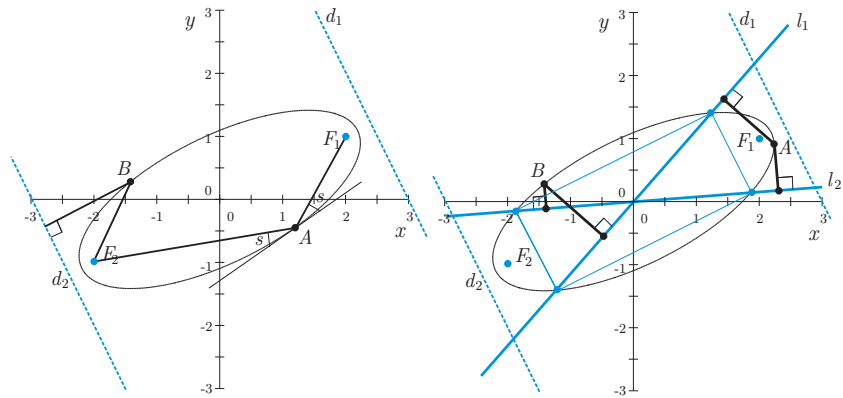}%
\caption{Two views of the ellipse $2x^{2}-4xy+5y^{2}=6$}%
\label{EllipseBlue}%
\end{center}
\end{figure}

The left diagram in Figure \ref{EllipseBlue} shows the central ellipse
\[
2x^{2}-4xy+5y^{2}=6
\]
with foci at $F_{1}=\left[  2,1\right]  $ and $F_{2}=\left[  -2,-1\right]  $,
corresponding directrices $d_{1}$ and $d_{2}$ with respective equations
$2x-y+6=0$ and $2x-y-6=0,$ and eccentricity $e=\sqrt{5/6}.$ The familiar
reflection property may be recast as: \textit{spreads} between a tangent and
lines to the foci from a point on the ellipse are equal.

The right diagram in Figure \ref{EllipseBlue} illustrates a (perhaps?) novel
definition of an ellipse. It is motivated by the fact that a circle is the
locus of a point $X$ whose quadrances to two fixed perpendicular lines add to
a constant. Define a \textbf{grammola} to be the locus of a point $X$ such
that the sum of the quadrances from $X$ to two fixed non-perpendicular
intersecting lines $l_{1}$ and $l_{2}$ is constant. This definition works for
each of the three colours. It turns out that the lines $l_{1}$ and $l_{2} $
are unique; we call them the \textbf{diagonals }of the grammola (see
\cite[Chapter 15]{Wild1}). The \textbf{corners} of the grammola are the points
where the diagonal lines intersect it, and determine the \textbf{corner
rectangle}. In the blue setting over the real numbers a grammola is always an
\textit{ellipse}, while in the red and green settings a grammola might be an
ellipse, or it might be a hyperbola.

The ellipse of Figure \ref{EllipseBlue} is a blue grammola with blue diagonals%
\[%
\begin{tabular}
[c]{lllll}%
$\left(  14+5\sqrt{6}\right)  x-23y=0$ &  & and &  & $\left(  14-5\sqrt
{6}\right)  x-23y=0.$%
\end{tabular}
\]
The blue quadrances of the sides of the corner rectangle are $12$ and $2,$
whose product $24$ is the squared area. The right diagram in Figure
\ref{EllipseBlue} shows the usual foci and directrices of the grammola and its
diagonals and corners. The quadrances from any point on the conic to the two
diagonals sum to $6.$ The blue spread between the two diagonals is an
invariant of the ellipse---in this case $s_{b}=24/49$.%
\begin{figure}
[h]
\begin{center}
\includegraphics[
height=5.162cm,
width=10.772cm
]%
{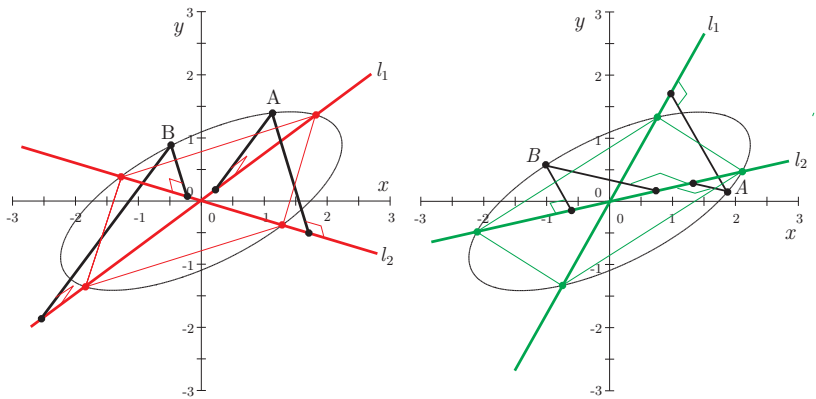}%
\caption{Ellipse as red and green grammola}%
\label{RedGreen}%
\end{center}
\end{figure}

The ellipse can also be described as a red grammola, as in the left diagram of
Figure \ref{RedGreen}. The red diagonals are%
\[%
\begin{tabular}
[c]{lllll}%
$\left(  \sqrt{22}+2\right)  x-9y=0$ &  & \textrm{and} &  & $\left(  \sqrt
{22}-2\right)  x+9y=0$%
\end{tabular}
\]
and the \textbf{red corner rectangle} has sides parallel to the red axes of
the ellipse, and red quadrances $3+\sqrt{33}$ and $3-\sqrt{33},$ whose product
is $-24$. The four red corners have rather complicated expressions in this
case. The red spread between the red diagonals is $s_{r}=-8/3.$

The same ellipse may also be viewed as a green grammola, as in the right
diagram of Figure \ref{RedGreen}. The green diagonals are
\[%
\begin{tabular}
[c]{lllll}%
$\left(  -5+\sqrt{15}\right)  x+5y=0$ &  & \textrm{and} &  & $\left(
-5-\sqrt{15}\right)  x+5y=0$%
\end{tabular}
\]
and the \textbf{green corner rectangle} has sides parallel to the green axes
of the ellipse, and green quadrances $4+2\sqrt{10}$ and $4-2\sqrt{10},$ whose
product is again $-24$. Except for a sign, the three squared areas of the
blue, red and green corner rectangles are the same. The green spread between
the green diagonals is $s_{g}=-3/2.$

The relationship between the blue, red and green spreads of an ellipse is
\[
\frac{1}{s_{b}}+\frac{1}{s_{r}}+\frac{1}{s_{g}}=1.
\]

\section*{The ellipse as a quadrola}

Another well known definition of an ellipse is as the locus of a point $X$
whose sum of distances from two fixed points $F_{1}$ and $F_{2}$ is a constant
$k$. To determine a universal analog of this, we consider the locus of a point
$X$ such that the quadrances $Q_{1}=Q\left(  F_{1},X\right)  $ and
$Q_{2}=Q\left(  F_{2},X\right)  $, together with a number $K,$ satisfy
Archimedes formula $A\left(  Q_{1},Q_{2},K\right)  =0.$ This is the quadratic
analog to the equation $d_{1}+d_{2}=k,$ just as the Triple quad formula is the
analog to a linear relation between three distances.

Such a locus we call a \textbf{quadrola}. This algebraic formulation applies
to the relativistic geometries, and also extends to general fields. The notion
captures both that of ellipse and hyperbola in the Euclidean setting, and
while it is in general a different concept than a grammola, it is possible for
a conic to be both, as is the case of an ellipse in Euclidean (blue) geometry.%
\begin{figure}
[h]
\begin{center}
\includegraphics[
height=5.1793cm,
width=10.9572cm
]%
{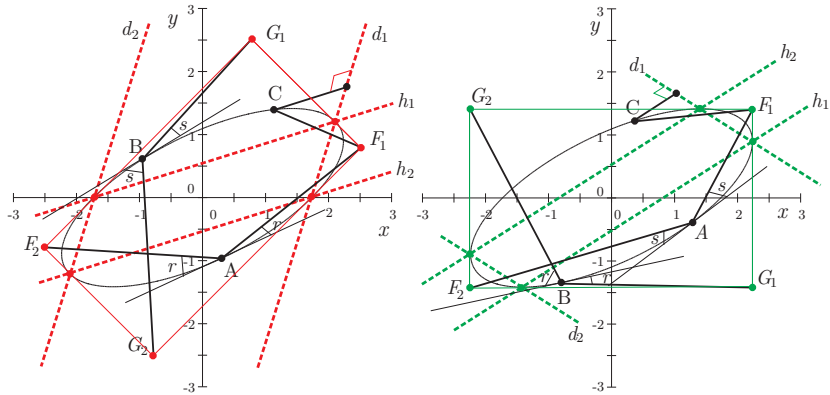}%
\caption{Red and green foci and directrices}%
\label{EllipseRedGreen}%
\end{center}
\end{figure}

The left diagram in Figure \ref{EllipseRedGreen} shows that in the red
geometry, a new phenomenon occurs: our same ellipse as a quadrola has
\textit{two pairs} of foci $\left\{  F_{1},F_{2}\right\}  $ and $\left\{
G_{1},G_{2}\right\}  $. Each of these points is also a focus in the context of
a conic section, and there are two pairs of corresponding directrices
$\left\{  d_{1},d_{2}\right\}  $ and $\left\{  h_{1},h_{2}\right\}  $.
Directrices are parallel or red perpendicular, and intersect at points on the
ellipse, and tangents to these \textbf{directrix points} pass through
two\textit{\ }foci, forming a parallelogram which is both a \textit{blue} and
a \textit{green} \textit{rectangle}. It turns out that the red spreads between
a tangent and lines to a pair of red foci are equal, as shown at points $A$
and $B$.

The right diagram in Figure \ref{EllipseRedGreen} shows the same ellipse as a
green quadrola, with again two pairs of green foci, two pairs of corresponding
green directrices (which are parallel or green perpendicular), and the
tangents at directrix points forming a blue and red rectangle.

The red and green directrix points are easy to find: they are the limits of
the ellipse in the null and the coordinate directions. So the red and green
directrices and foci are also then simple to locate geometrically. This is not
the case for the usual (blue) foci and directrices, and suggests that
considering ellipses from the relativistic perspectives can be practically
useful. In algebraic geometry the `other' pair of blue foci are not unknown;
they require complexification and a projective view (see for example
\cite[Chapter 12]{Gibson}).%
\begin{figure}
[h]
\begin{center}
\includegraphics[
height=5.1944cm,
width=5.4678cm
]%
{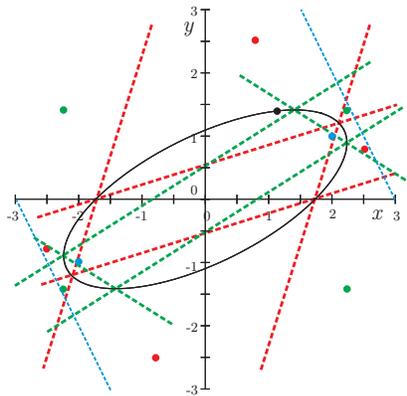}%
\caption{Three sets of foci and directrices}%
\label{ThreeQuadrola}%
\end{center}
\end{figure}

When we put all three coloured pictures together, as shown in Figure
\ref{ThreeQuadrola}, another curious phenomenon appears---there are three
pairs of coloured foci that appear to be close to the intersections of
directrices of the opposite colour. The reason for this will become clearer
later when we consider parabolas.

\section*{Hyperbolas}

\begin{figure}
[h]
\begin{center}
\includegraphics[
height=5.2762cm,
width=5.2417cm
]%
{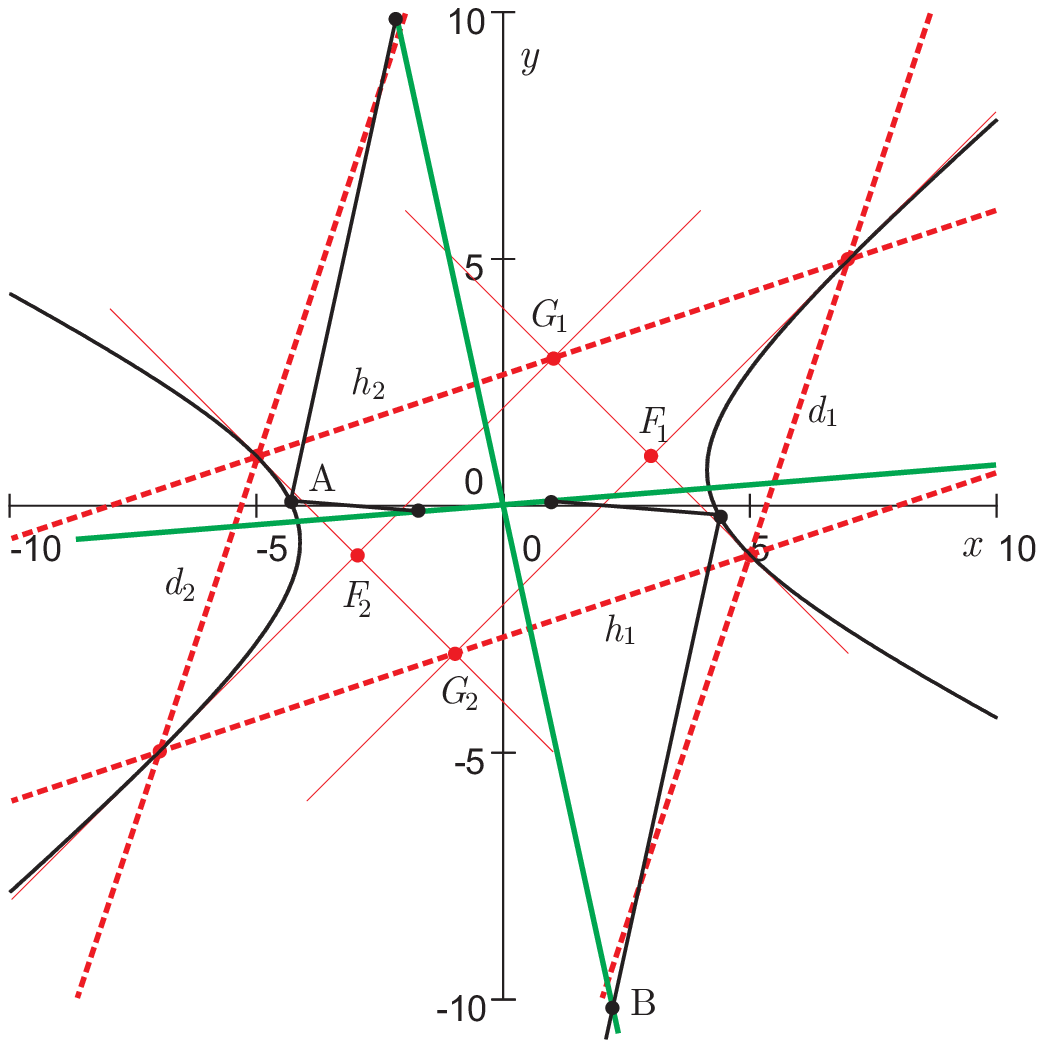}%
\caption{A hyperbola as red quadrola and green grammola}%
\label{RedHyperbola}%
\end{center}
\end{figure}

Over the real numbers some of what we saw with ellipses extends also to
hyperbolas, although there are differences. The central hyperbola shown in
Figure \ref{RedHyperbola} with equation
\[
7x^{2}+6xy-17y^{2}=128
\]
is a red quadrola with red foci $F_{1}=\left[  3,1\right]  $ and
$F_{2}=\left[  -3,-1\right]  $, meaning that it is the locus of a point
$X=\left[  x,y\right]  $ such that
\[
A\left(  \left(  x-3\right)  ^{2}-\left(  y-1\right)  ^{2},\left(  x+3\right)
^{2}-\left(  y+1\right)  ^{2},64\right)  =0.
\]%

As a conic section the corresponding directrices have equations $3x-y-16=0$
and $3x-y+16=0.$ This hyperbola also has another pair of red foci
$G_{1}=\left[  1,3\right]  $ and $G_{2}=\left[  -1,-3\right]  $, with
associated directrices $x-3y+8=0$ and $x-3y-8=0.$

As in the case of the ellipse we considered earlier, in each case the focus is
the pole of the corresponding directrix, meaning that it is the intersection
of the tangents to the hyperbola at the directrix points. These tangents pass
through two foci at a time, and are parallel to the red null directions. The
parallelogram formed by the four foci is a blue and green rectangle.

So we could have found the red foci and directrices purely geometrically, by
finding those points on the hyperbola where the tangents are parallel to the
red null directions, and then forming intersections between these points. This
is again quite different from finding the usual blue foci and directrices.
Note that if we try to find green foci, vertical tangents are easy to find,
but there are no horizontal tangents, thus the situation will necessarily be
somewhat different.

Is the hyperbola also a grammola? It cannot be a blue grammola, since these
are all ellipses, and it turns out not to be a red grammola either. But it
\textit{is} a green grammola with equation
\[
\frac{\left(  \left(  119+8\sqrt{238}\right)  x+51y\right)  ^{2}}{2\left(
119+8\sqrt{238}\right)  51}+\frac{\left(  \left(  119-8\sqrt{238}\right)
x+51y\right)  ^{2}}{2\left(  119-8\sqrt{238}\right)  51}=\frac{128}{3}.
\]
The green diagonals are shown in Figure \ref{RedHyperbola}.

\section*{The parabola}

From the viewpoint of universal (affine) geometry, the most interesting conic
is the \textit{parabola}. Given a point $F$ and a generic line $l$ not passing
through $F$, the locus of a point $X$ such that $Q\left(  X,F\right)
=Q\left(  X,l\right)  $ is what we usually call a parabola, independent of
which geometry we are considering. The generic parabola has a distinguished
blue, red and green focus, and also a blue, red and green directrix.
\begin{figure}
[h]
\begin{center}
\includegraphics[
height=5.3925cm,
width=6.8569cm
]%
{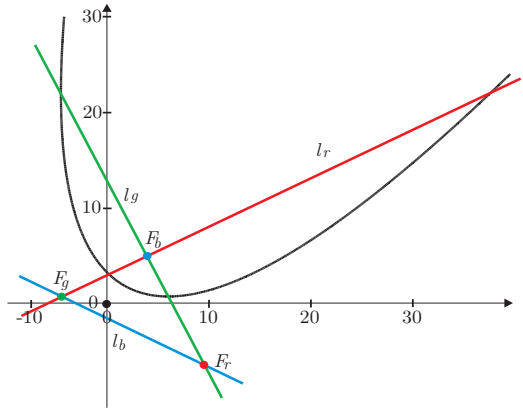}%
\caption{Three foci and directrices of a parabola}%
\label{ParabolaFoci}%
\end{center}
\end{figure}

Figure \ref{ParabolaFoci} shows a parabola in the Cartesian plane and all
three foci and directrices. A remarkable phenomenon appears: $F_{b}$ is the
intersection of $l_{r}$ and $l_{g}$, $F_{r}$ is the intersection of $l_{b}$
and $l_{g},$ and $F_{g}$ is the intersection of $l_{b}$ and $l_{r}$.
Furthermore $l_{r}$ and $l_{g}$ are blue perpendicular, $l_{b}$ and $l_{g}$
are red perpendicular, and $l_{b}$ and $l_{r}$ are green perpendicular---in
other words we get a triple right triangle of foci. This means that once we
know one of the focus/directrix pairs, the other two can be found simply by
constructing the appropriate altitudes from the focus to the directrix
together with their feet.
\begin{figure}
[h]
\begin{center}
\includegraphics[
height=5.3946cm,
width=6.8569cm
]%
{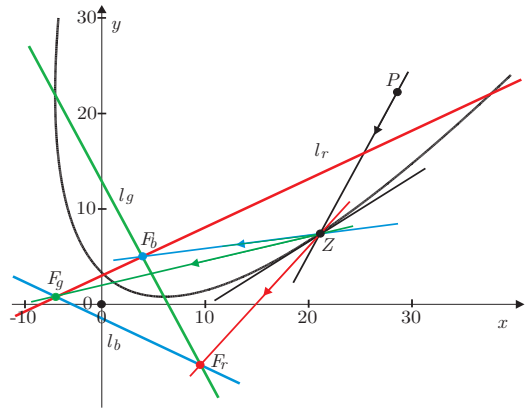}%
\caption{Reflection properties of a parabola}%
\label{Reflection}%
\end{center}
\end{figure}

Although the various directrices are in different directions, the
\textit{axis} direction, defined as being perpendicular to the directrix, is
common to all. Figure \ref{Reflection} shows the familiar reflection property
of the parabola, where a particle $P$ approaching the parabola along the axis
direction and reflecting off the tangent (in either a blue, red or green
fashion) always then passes through the corresponding focus.

The following figure shows some interesting collinearities associated to a
parabola, involving coloured vertices $V$ (intersections of axes with the
parabola), bases $X$ (intersections of axes with directrices) and points $Y$
formed by tangents to vertices.
\begin{figure}
[h]
\begin{center}
\includegraphics[
height=6.4456cm,
width=8.2524cm
]%
{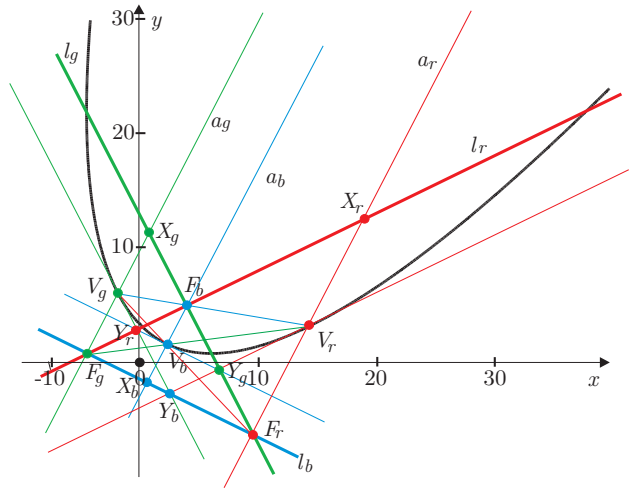}%
\caption{Collinearities for a chromatic parabola}%
\label{Vertices}%
\end{center}
\end{figure}

Finally we show the three parabolas which have a given focus $F$ and a given
directrix $l,$ both in black, each interpreted in one of the three geometries.
Each of the three parabolas that share this focus and directrix have a focal
triangle consisting of $F$ and two of the feet of the altitudes from $F$ to
$l,$ labelled $F_{b},F_{r}$ and $F_{g}.$ The dotted line passes through the
intersections of the red and green parabolas. Various vertices and axes are
shown, and we leave the reader to notice interesting collinearities, and to
try to prove them.%

\begin{figure}
[h]
\begin{center}
\includegraphics[
height=6.7578cm,
width=8.7218cm
]%
{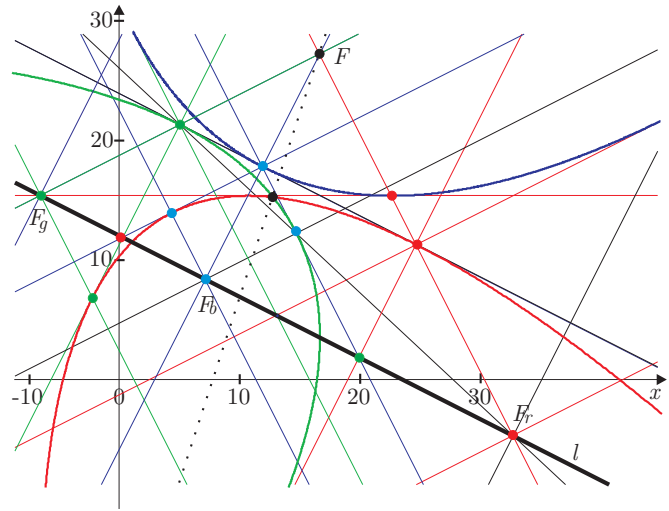}%
\caption{Three parabolas with a common focus and directrix}%
\label{ThreeParabolas}%
\end{center}
\end{figure}

In conclusion, there may very well be other useful metrical definitions of
conics; there are certainly still many rich discoveries to be made about these
fascinating and most important geometric objects. Chromogeometry extends to
many other aspects of planar geometry, for example to triangle geometry in
\cite{Wild3}.

\end{document}